\documentclass[11pt, a4paper]{amsart}
\usepackage{amsmath, amsthm, amssymb}
\usepackage{txfonts, mathrsfs}
\usepackage[usenames]{color}
\usepackage{psfrag}
\usepackage{graphicx}

\numberwithin{equation}{section}

\newtheorem{thm}{Theorem}[section]
\newcommand{\bt}{\begin{thm}}
\newcommand{\et}{\end{thm}}

\newtheorem{cor}[thm]{Corollary}
\newcommand{\bc}{\begin{cor}}
\newcommand{\ec}{\end{cor}}

\newtheorem{lem}[thm]{Lemma}
\newcommand{\bl}{\begin{lem}}
\newcommand{\el}{\end{lem}}

\newtheorem{prop}[thm]{Proposition}
\newcommand{\bp}{\begin{prop}}
\newcommand{\ep}{\end{prop}}

\newtheorem{defn}[thm]{Definition}
\newcommand{\bd}{\begin{defn}}
\newcommand{\ed}{\end{defn}}

\newtheorem{rmrk}[thm]{Remark}
\newcommand{\br}{\begin{rmrk}}
\newcommand{\er}{\end{rmrk}}

\newtheorem{quest}[thm]{Question}
\newcommand{\bq}{\begin{quest}}
\newcommand{\eq}{\end{quest}}

\newcommand{\C}{\mathbb{C}}
\newcommand{\N}{\mathbb{N}}

\newdimen\vintkern\vintkern12pt
\def\vint{-\kern-\vintkern\int}

\newcommand{\hm}{{\mathcal H}}

\newcommand{\trace}{\operatorname{tr}}

\newcommand{\Area}{\operatorname{Area}}

\newcommand{\md}{\operatorname{md}}

\newcommand{\jac}{{\mathbf J}}
\newcommand{\ap}{\operatorname{ap}}
\newcommand{\apmd}{\ap\md}

\newcommand{\R}{\mathbb{R}}

\begin{document}
\pagebreak
\bibliographystyle{plain}


\title[Morrey's $\varepsilon$--conformality lemma]{Morrey's $\varepsilon$--conformality lemma in metric spaces}


\author{Martin Fitzi}

\address
  {Department of Mathematics\\ University of Fribourg\\ Chemin du Mus\'ee 23\\ 1700 Fribourg, Switzerland}
\email{martin.fitzi@unifr.ch}

\author{Stefan Wenger}

\address
  {Department of Mathematics\\ University of Fribourg\\ Chemin du Mus\'ee 23\\ 1700 Fribourg, Switzerland}
\email{stefan.wenger@unifr.ch}

\date{\today}

\thanks{This work was partially supported by the following grants: Swiss National Science Foundation Grants 165848 and 182423;  Grant 346300 for IMPAN from the Simons Foundation and the matching 2015--2019 Polish MNiSW fund.}

\begin{abstract}
We provide a simpler proof and slight strengthening of Morrey's famous lemma on $\varepsilon$--conformal mappings. Our result more generally applies to Sobolev maps with values in a complete metric space and we obtain applications to the existence of area minimizing surfaces of higher genus in metric spaces. Unlike Morrey's proof, which relies on the measurable Riemann mapping theorem, we only need the existence of smooth isothermal coordinates established by Korn and Lichtenstein.
\end{abstract}

\maketitle

\renewcommand{\theequation}{\arabic{section}.\arabic{equation}}
\pagenumbering{arabic}

\section{Introduction and statement of results}\label{sec:Intro}

\subsection{Introduction}\label{sec:intro-basics}

Let $N$ be a Riemannian manifold, denote by $D$ the open unit disc in $\R^2$, and let $u\in W^{1,2}(D, N)$ be a Sobolev map. The Dirichlet energy of $u$ and the parametrized area, i.e.~the integral of the jacobian determinant, satisfy the inequality
\begin{equation}\label{eq:area-energy-Rn}
\Area(u) \leq E^2(u).
\end{equation}
Furthermore, equality holds if and only if $u$ is weakly conformal in the sense that almost everywhere the weak partial derivatives $u_x$ and $u_y$ are orthogonal and have the same length (possibly zero).
An important result, proved by Morrey in \cite{Mor48} and known as Morrey's lemma on $\varepsilon$--conformal mappings, shows that after a suitable reparametrization of $u$ one can obtain near equality in \eqref{eq:area-energy-Rn}. More precisely, for every $\varepsilon>0$ there exists a homeomorphism $\varphi$ of the disc $D$ such that $u\circ\varphi\in W^{1,2}(D, N)$ and $$E^2(u\circ\varphi)\leq \Area(u) + \varepsilon.$$ The proof relies on the existence of $L^2$--solutions of the Beltrami equation with measurable coefficients, previously established by Morrey in \cite{Mor38} and called the measurable Riemann mapping theorem. The homeomorphism $\varphi$ is actually a quasiconformal map in the sense of geometric function theory.

In this note we prove a version of Morrey's lemma for Sobolev maps defined on a smooth, compact $2$-dimensional manifold and with values in a metric space. Our proof is similar to that of Morrey but contains some simplifications despite the greater generality. Specifically, we do not use the measurable Riemann mapping theorem but only need the existence of isothermal coordinates established by Korn \cite{Kor19} and Lichtenstein \cite{Lic16}, thus solutions to the Beltrami equation with smooth compactly supported coefficients. We actually obtain a stronger statement than asserted by the classical lemma of Morrey. Our result has applications to the existence of area minimizing surfaces of higher genus and the structure of energy minimizing surfaces.

\subsection{Main results} 

Throughout this introduction, let $X$ be a complete metric space and let $M$ be a smooth, compact, connected, orientable $2$-dimensional manifold, possibly with boundary. 
We denote by $W^{1,2}(M,X)$ the space of $2$--Sobolev maps from $M$ to $X$ in the sense of Reshetnyak \cite{Res97}, \cite{Res06}. There exist several different but equivalent definitions of Sobolev maps from Euclidean or Riemannian domains into a metric space. For references as well as for definitions of the following notions see Section~\ref{sec:Sobolev-prelims} below. 
Given a Riemannian metric $g$ on $M$ we denote by $E_+^2(u,g)$ the Reshetnyak energy of $u\in W^{1,2}(M, X)$ with respect to $g$. We furthermore let $\Area_{\mu^i}(u)$ be the intrinsic Riemannian area  of $u$. This definition of area was first studied by Ivanov \cite{Iva08} and later in \cite{LW17-en-area}; it agrees with the usual parametrized area when $X$ is a Riemannian manifold. In the setting of non-Euclidean metric spaces there are several natural definitions of energy and area. The strong implication in the equality case in the following proposition established in \cite{LW15-Plateau}, \cite{LW17-en-area} is one of the reasons for our choice of energy and area.

\bp\label{prop:energy-area-ineq-inf-isotrop}
For every $u\in W^{1,2}(M, X)$ and every Riemannian metric $g$ on $M$ we have $$\Area_{\mu^i}(u)\leq E_+^2(u,g),$$ with equality if and only if $u$ is infinitesimally isotropic with respect to $g$.
\ep

The precise meaning of {\it infinitesimally isotropic} will be explained in Definition~\ref{def:inf-isotropic}. This notion provides a substitute for weak conformality in the setting of metric spaces and it is equivalent to weak conformality when $X$ is a Riemannian manifold. For general metric spaces, if $u$ is infinitesimally isotropic then it is in particular infinitesimally $\sqrt{2}$--quasiconformal. That is, up to a factor of $\sqrt{2}$, the infinitesimal metric distortion under $u$ of all unit vectors at a given point is the same; see \eqref{eq:def-inf-sqrt2-qc} for the precise definition.

The following is one of our main results:

\bt\label{thm:Morrey-epsilon-conformality}
 For every $u\in W^{1,2}(M, X)$ and every $\varepsilon>0$ there exists a Riemannian metric $g$ on $M$ such that $$E_+^2(u, g) \leq \Area_{\mu_i}(u) + \varepsilon.$$
\et

The metric $g$ can be chosen to have constant sectional curvature $-1$, $0$, or $1$ and such that the boundary of $M$ is geodesic. Theorem~\ref{thm:Morrey-epsilon-conformality} and Proposition~\ref{prop:energy-area-ineq-inf-isotrop} yield the following corollary which strengthens \cite[Theorem 4.2]{FW-Plateau-Douglas}.

\bc\label{cor:energy-min-inf-isotropic}
 Let $u\in W^{1,2}(M, X)$ and let $g$ be a Riemannian metric on $M$. If for every Riemannian metric $h$ on $M$ we have $$E_+^2(u,g)\leq E_+^2(u,h)$$ then $u$ is infinitesimally isotropic with respect to $g$.
\ec

The proof of Theorem~\ref{thm:Morrey-epsilon-conformality} relies on the following analog for Sobolev maps defined on the open unit disc $D\subset\R^2$ and on the uniformization theorem for compact Riemann surfaces. For a Sobolev map $u\in W^{1,2}(D,X)$ we write $E_+^2(u)$ for its Reshetnyak energy with respect to the Euclidean metric on $D$.

\bt\label{thm:Morrey-epsilon-Reshetnyak-disc-intro}
For every $u\in W^{1,2}(D,X)$ and every $\varepsilon>0$ there exists a diffeomorphism $\varphi\colon D\to D$ such that $$E_+^2(u\circ\varphi) \leq \Area_{\mu_i}(u) + \varepsilon.$$
Moreover, the map $\varphi$ extends to a diffeomorphism of the closed unit disc $\overline{D}$ and is conformal in a neighborhood of the boundary.
\et

When $X$ is a Riemannian manifold then the Reshetnyak energy bounds from above the Dirichlet energy, see Section~\ref{sec:Sobolev-prelims}, so Theorem~\ref{thm:Morrey-epsilon-Reshetnyak-disc-intro} in particular implies and strengthens the classical lemma of Morrey \cite[Theorem 1.2]{Mor48}. Together with Proposition~\ref{prop:energy-area-ineq-inf-isotrop} we furthermore obtain the following result which implies \cite[Theorem 6.2]{LW15-Plateau}.

\bc\label{bc:energy-min-inf-isotropic-disc}
 If $u\in W^{1,2}(D,X)$ satisfies $E_+^2(u)\leq E_+^2(u\circ\varphi)$ for every diffeomorphism $\varphi$ of the closed unit disc $\overline{D}$ then $u$ is infinitesimally isotropic (with respect to the Euclidean metric).
\ec

Our results are useful in the context of the Plateau-Douglas problem. This problem is concerned with the existence of area minimizing surfaces of fixed genus spanning a given finite collection of Jordan curves. In the setting of Euclidean space or, more generally, Riemannian manifolds the problem has been solved under two different conditions of non-degeneracy: the so-called Douglas condition considered in \cite{Dou39}, \cite{Jos85}, as well as Courant's condition of cohesion used in \cite{Shi39}, \cite{Cou40}, \cite{TT88}. Assuming the Douglas condition we have recently solved the Plateau-Douglas problem in any proper metric space admitting a local quadratic isoperimetric inequality in \cite{FW-Plateau-Douglas}. In the present paper we solve this problem using Courant's condition of cohesion in any proper metric space.

Let $\Gamma\subset X$ be the disjoint union of $k\geq 1$ rectifiable Jordan curves in $X$ and suppose the smooth surface $M$ has $k$ boundary components. We denote by $\Lambda(M,\Gamma,X)$ the possibly empty family of maps $u\in W^{1,2}(M, X)$ whose trace has a continuous representative which weakly monotonically parametrizes $\Gamma$, i.e.~is the uniform limit of homeomorphisms $\varphi_j\colon \partial M\to \Gamma$. We can now state our solution to the Plateau-Douglas problem, which generalizes the results in \cite{Shi39}, \cite{Cou40}, \cite{TT88} to the setting of proper metric spaces.

\bt\label{thm:plateau-douglas-cond-cohesion}
 Suppose the metric space $X$ is proper, $\Lambda(M, \Gamma, X)$ is not empty and contains a $\mu^i$--area minimizing sequence which satisfies the condition of cohesion. Then there exist $u\in\Lambda(M, \Gamma, X)$ and a Riemannian metric $g$ on $M$ such that $$\Area_{\mu^i}(u) =  \inf\{\Area_{\mu^i}(v): v\in\Lambda(M, \Gamma, X)\}$$ and such that $u$ is infinitesimally isotropic with respect to $g$.
\et

We refer to Section~\ref{sec:area-min-cond-cohesion} for the definition of Courant's condition of cohesion. The Riemannian metric $g$ can be chosen in such a way that $(M,g)$ has constant curvature $-1$, $0$ or $1$ and that $M$ has geodesic boundary. If the metric space $X$ admits a local quadratic isoperimetric inequality for curves then the map $u$ in Theorem~\ref{thm:plateau-douglas-cond-cohesion} has a unique representative which is locally H\"older continuous in the interior of $M$ and is continuous up to the boundary, see \cite[Theorem 1.4]{FW-Plateau-Douglas}.

Recall that in the generality of metric spaces, there are several natural notions of parametrized area of a Sobolev map. It is natural to ask whether Theorem~\ref{thm:plateau-douglas-cond-cohesion} holds when the intrinsic Riemannian area is replaced by another notion of area, for example the parametrized Hausdorff area. In Section~\ref{sec:area-min-cond-cohesion} we will show that there exists a Hausdorff area minimizer in $\Lambda(M, \Gamma, X)$ but at present we do not know how to get a minimizer with a good parametrization as in the theorem above.

\bigskip

{\bf Acknowledgements:} We thank Katrin F\"assler, Pekka Pankka, Kai Rajala and Teri Soultanis for useful comments.

\section{Preliminaries}

\subsection{Basic notation}

We write $|v|$ for the Euclidean norm of a vector $v\in\R^2$, $$D:=\{z\in\R^2:|z|<1\}$$ for the open unit disc in $\R^2$ and $\overline{D}$ for its closure.
The differential at $z$ of a differentiable map $\varphi$ between smooth manifolds is denoted $D\varphi(z)$.

For a subset $A\subset\R^2$, $|A|$ denotes its Lebesgue measure. If $(X,d)$ is a metric space then we use the notation $\hm_{X}^2(A)$ for the $2$--dimensional Hausdorff measure of a subset $A\subset X$. The normalizing constant is chosen such that $\hm_X^2$ coincides with the $2$--dimensional Lebesgue measure when $X$ is Euclidean $\R^2$. Thus, the Hausdorff $2$--measure $\hm_g^2:=\hm_{(M,g)}^2$ on a 2--dimensional Riemannian manifold $(M,g)$ coincides with the Riemannian area.

\subsection{Energy, area, and isotropy of semi-norms}\label{sec-energy-area-isotropy}
The (Reshetnyak) energy of a semi-norm $s$ on $\R^2$ is defined by  $$\mathbf{I}_+^2(s):= \max\{s(v)^2: v\in\R^2, |v|=1\}.$$
Next, we define the jacobian $\jac_{\mu^i}$ corresponding to the intrinsic Riemannian or Loewner area. For a degenerate semi-norm $s$ we set $\jac_{\mu^i}(s):=0$. If $s$ is a norm on $\R^2$ and $L$ is the ellipse of maximal area contained in $\{v\in\R^2: s(v)\leq 1\}$ we set $$\jac_{\mu^i}(s) = \frac{\pi}{|L|},$$ where $|L|$ is the Lebesgue measure of $L$. Thus, the measure $$\mu^i_{(\R^2, s)}(A):= \jac_{\mu^i}(s)\cdot |A|$$ is the multiple of the Lebesgue measure for which $L$ has measure $\pi$. In particular, if $s$ is induced by an inner product then $\mu^i_{(\R^2, s)}(A)$ agrees with the Hausdorff $2$--measure $\hm^2_{(\R^2, s)}(A)$.

\bd\label{def:isotropic-seminorm}
 A semi-norm $s$ on $\R^2$ is isotropic if $s=0$ or if it is a norm and the ellipse of maximal area contained in $\{v\in\R^2: s(v)\leq 1\}$ is a Euclidean ball.
\ed

It follows from \cite[Section 3.2]{LW17-en-area} that 
\begin{equation}\label{eq:jac-energy}
\jac_{\mu^i}(s) \leq \mathbf{I}_+^2(s)
\end{equation}
 for every semi-norm on $\R^2$, with equality if and only if $s$ is isotropic.

\subsection{Approximate metric derivatives}\label{sec:aprox-metr-deriv}

A map $v\colon V\to X$ from an open subset $V\subset \R^2$ to a metric space $(X,d)$ is said to be approximately metrically differentiable at a point $z\in V$ if there exists a semi-norm $s$ on $\R^2$ such that $$\ap\lim_{z'\to z} \frac{d(v(z'), v(z)) - s(z'-z)}{|z'-z|} = 0,$$ where $\ap\lim$ denotes the approximate limit, see for example \cite{EG92}. If such a semi-norm exists then it is unique and is called the approximate metric derivative of $v$ at $z$ and denoted $\apmd v_z$. 

Let $M$ be a $2$--dimensional smooth manifold without boundary and $z\in M$. A map $u\colon M\to X$ is called approximately metrically differentiable at $z$ if the composition $u\circ\psi^{-1}$ is approximately metrically differentiable at $\psi(z)$ for some and thus every chart $(U, \psi)$ of $M$ around $z$. The semi-norm on the tangent space $T_zM$ defined by $$\apmd u_z:= \apmd (u\circ\psi^{-1})_{\psi(z)}\circ D\psi(z)$$ is independent of the choice of chart and is called the approximate metric derivative of $u$ at $z$.

Let $g$ be a Riemannian metric on $M$ and fix $z\in M$. If $s$ is a semi-norm on $V:= (T_zM, g(z))$ then we define the concepts of jacobian, energy, and isotropy introduced in Section~\ref{sec-energy-area-isotropy} by identifying $V$ with Euclidean $(\R^2,|\cdot|)$ via a linear isometry.

\subsection{Beltrami coefficients and Beltrami equation}\label{sec:Beltrami}

In the proof of Theorem~\ref{thm:Morrey-epsilon-Reshetnyak-disc-intro} we will need diffeomorphic solutions to the Beltrami equation $f_{\overline{z}} = \mu f_z$ for smooth compactly supported Beltrami coefficients $\mu$. For such coefficients, the solution goes back to work of Korn \cite{Kor19} and Lichtenstein \cite{Lic16} on the existence of isothermal coordinates. In this section we recall the Beltrami equation and the main existence theorem which we will use.  For details we refer for example to \cite{AB60} or \cite{AIM09}. It will be convenient to work with complex notation.

Let $f\colon\C\to\C$ be an orientation preserving diffeomorphism. We use the complex differential operators $$f_z = \frac{1}{2}(f_x - i f_y) \quad\text{ and } \quad f_{\overline{z}} = \frac{1}{2}(f_x + i f_y),$$ where $f_x$ and $f_y$ are the partial derivatives of $f$ with respect to $x$ and $y$.
The differential $Df(z)$ of $f$ at a point $z\in\C$ can be written as $$Df(z)(h) = f_z(z)\cdot h + f_{\overline{z}}(z)\cdot \overline{h}$$ for every $h\in\C$. It is not difficult to show that the operator norm of the differential satisfies
$$\|Df(z)\| = \max_{|h|=1} |Df(z)(h)| =  |f_z(z)| + |f_{\overline{z}}(z)|.$$ We moreover have $$\det(Df(z)) = |f_z(z)|^2 - |f_{\overline{z}}(z)|^2>0$$ as well as $$\min_{|h|=1} |Df(z)(h)| = |f_z(z)| - |f_{\overline{z}}(z)|.$$
We thus obtain the identities
\begin{equation}\label{eq:Beltrami-to-distortion}
 \frac{\|Df(z)\|^2}{\det(Df(z))} = \frac{\max_{|h|=1} |Df(z)(h)|}{\min_{|h|=1} |Df(z)(h)|} = \frac{1+|\mu_f(z)|}{1-|\mu_f(z)|},
\end{equation}
where $\mu_f(z)$ is the {\it Beltrami coefficient} of $f$ at $z$ defined by the equation $$f_{\overline{z}}(z)=\mu_f(z)f_z(z).$$
Given $K\geq 1$, the orientation preserving diffeomorphism $f$ is said to be $K$--quasi\-conformal if $$\|Df(z)\|^2 \leq K\cdot \det(Df(z))$$ for every $z\in\C$.
If follows from \eqref{eq:Beltrami-to-distortion} that $f$ is $K$--quasiconformal if and only if the $L^\infty$--norm of the Beltrami coefficient satisfies $$\|\mu_f\|_\infty\leq \frac{K-1}{K+1}<1.$$ Moreover, if $f$ is $K$--quasiconformal then its inverse $f^{-1}$ is also $K$--quasiconformal.

The following existence result goes back to Korn \cite{Kor19} and Lichtenstein \cite{Lic16}, see the bibliographical note in \cite{AB60}. For a proof we refer for example to \cite{AB60} or \cite[Theorem 5.2.4]{AIM09}.

\bt\label{thm:Sol-Beltrami-qc-diffeo}
 Let $\mu\colon \C\to\C$ be a smooth and compactly supported function satisfying $k:= \|\mu\|_\infty<1$. Then there exists an orientation preserving diffeomorphism $f\colon \C\to\C$ which is $K$--quasiconformal with $K= \frac{1+k}{1-k}$ and which solves the Beltrami equation $f_{\overline{z}} = \mu f_z$ in $\C$.
\et

We will need the following composition formula for the Beltrami coefficient. Let $f,g\colon \C\to \C$ be two orientation preserving diffeomorphisms. Then the Beltrami coefficient of the composition $g\circ f^{-1}$ at a point $w=f(z)$ is given by 
\begin{equation}\label{eq:composition-Beltrami-coeff}
 \mu_{g\circ f^{-1}}(w) = \frac{\mu_g(z) - \mu_f(z)}{1- \mu_g(z)\overline{\mu_f(z)}}\cdot \left(\frac{f_z(z)}{|f_z(z)|}\right)^2,
\end{equation}
see \cite[Theorem 5.5.6]{AIM09}.
 Moreover, if $f = \rho\circ g$ for some conformal diffeomorphism then $\mu_f \equiv \mu_g$.

\section{Sobolev maps with values in a metric space}\label{sec:Sobolev-prelims}

One can define Sobolev maps from a Euclidean or Riemannian domain into a complete metric space in several equivalent ways, see for example \cite{Amb90}, \cite{KS93}, \cite{Haj96}, \cite{Res97}, \cite{Jos97}, \cite{HKST15}. In this paper, we work with the definition introduced by Reshetnyak \cite{Res97}, \cite{Res06}, which we review here. We furthermore recall the definition of the Reshetnyak energy, the intrinsic Riemannian area and the notion of infinitesimal isotropy of a Sobolev map.

Let $(X, d)$ be a complete metric space and let $M$ be a smooth, compact, connected, orientable $2$--dimen\-sional manifold, possibly with non-empty boundary.
We fix a Riemannian metric $g$ on $M$ and let  $\Omega\subset M$ be an open set. 

\bd
A measurable and essentially separably valued map $u\colon \Omega\to X$ belongs to the Sobolev space $W^{1,2}(\Omega, X)$ if there exists $h\in L^2(\Omega)$ with the following property. For every real-valued $1$--Lipschitz function $f$ on $X$ the composition $f\circ u$ belongs to the classical Sobolev space $W^{1,2}(\Omega\setminus \partial M)$  and $$|\nabla (f\circ u)|\leq h$$ almost everywhere on $\Omega$. Here, $|\nabla (f\circ u)|$ denotes the length of the weak gradient of $f\circ u$ with respect to the metric $g$.
\ed

By \cite[Proposition 4.3]{LW15-Plateau}, a map $u\in W^{1,2}(\Omega, X)$ is approximately metrically differentiable at almost every $z\in M$ in the sense defined in Section~\ref{sec:aprox-metr-deriv} above. This allows one to define the Reshetnyak energy and the intrinsic Riemannian area using the pointwise quantities introduced in Section~\ref{sec-energy-area-isotropy} above.

\bd
 The Reshetnyak energy of $u\in W^{1,2}(\Omega, X)$ with respect to $g$ is defined by $$E_+^2(u, g):= \int_{\Omega} \mathbf{I}_+^2(\apmd u_z)\,d\hm^2_{g}(z).$$
\ed

The energy of the restriction of $u$ to a measurable subset $A\subset \Omega$ is defined in an analogous way. If $(U, \psi)$ is a conformal chart of $M$ with respect to $g$ then it follows from the area formula that $$E_+^2(u|_K, g) = \int_{\psi(K)} \mathbf{I}_+^2(\apmd(u\circ\psi^{-1})_w)\,dw = E_+^2(u\circ\psi^{-1}|_{\psi(K)}, g_{\rm Eucl})$$ for every compact set $K\subset U$, where $g_{\rm Eucl}$ denotes the Euclidean metric. This implies that the energy $E_+^2$ is invariant under precompositions with conformal diffeomorphisms.

It follows from \cite[Theorem 7.1.20]{HKST15} that $E_+^2(u,g)$ equals the square of the $L^2$--norm of the minimal weak upper gradient of a suitable representative of $u$. If $X$ is a Riemannian manifold then a short calculation shows that $E_+^2(u, g)$ bounds from above the classical Dirichlet energy $E^2(u, g)$ of $u$. Recall that if $\Omega$ is a subset of Euclidean $\R^2$ then $$E^2(u) = \frac{1}{2}\int_\Omega |u_x|^2 + |u_y|^2,$$ where $|u_x|$ and $|u_y|$ are the lengths of the weak partial derivatives of $u$ in the Riemannian manifold $X$. Using conformal charts and the conformal invariance of $E^2$ this definition easily extends to open subsets $\Omega$ of the $2$--dimensional smooth Riemannian manifold $(M,g)$.

\bd
 The intrinsic Riemannian area of $u\in W^{1,2}(\Omega, X)$ is defined by 
 \begin{equation}\label{eq:def-area}
  \Area_{\mu^i}(u):= \int_{\Omega} \jac_{\mu^i}(\apmd u_z)\,d\hm^2_g(z).
 \end{equation}
\ed

If $(U, \psi)$ is any chart of $M$ and $K\subset U$ is compact then $$\Area_{\mu^i}(u|_K) = \int_{\psi(K)} \jac_{\mu^i}(\apmd (u\circ\psi^{-1})_w)\,dw = \Area_{\mu^i}(u\circ\psi^{-1}|_{\psi(K)})$$ by the area formula. As a consequence, the parametrized $\mu^i$--area of a Sobolev map is invariant under precompositions with biLipschitz homeo\-morphisms. Finally, if $X$ is a Riemannian manifold or, more generally, a metric space with property (ET) of `Euclidean tangents' in the sense of \cite{LW15-Plateau} then the $\mu^i$--area agrees with the parametrized Hausdorff area, see Section~\ref{sec:area-min-cond-cohesion} below.

The following definition appears in \cite{LW17-en-area} and implicitly in \cite{LW15-Plateau} in the case that $(M,g)$ is a two-dimensional Euclidean domain.

\bd\label{def:inf-isotropic}
 A map $u\in W^{1,2}(M, X)$ is infinitesimally isotropic with respect to a Riemannian metric $g$ on $M$ if for almost every $z\in M$ the approximate metric derivative $\apmd u_z$ is isotropic.
\ed

If $u$ is infinitesimally isotropic with respect to $g$ then $u$ is infinitesimally $\sqrt{2}$--quasiconformal with respect to $g$ in the following sense. For almost all $z\in M$ we have 
\begin{equation}\label{eq:def-inf-sqrt2-qc}
\apmd u_z(v) \leq \sqrt{2}\cdot \apmd u_z(w)
\end{equation}
 for all $v, w\in T_zM$ with $|v|_z=|w|_z$ where $|v|_z$ is the length of $v$ with respect to $g$. If $X$ is a Riemannian manifold or, more generally, if $X$ has property (ET) of `Euclidean tangents' in the sense of \cite{LW15-Plateau} then the factor $\sqrt{2}$ can be replaced by $1$. 

\br\label{rem:area-engery-equality-inf-isotrop}
It follows from \eqref{eq:jac-energy} that 
\begin{equation}\label{eq:area-energy-ineq}
 \Area_{\mu^i}(u) \leq E_+^2(u,g)
\end{equation}
 for all Riemannian metrics $g$, and equality holds if and only if $u$ is infinitesimally isotropic.
\er

Finally, we recall the definition of the trace of a Sobolev map. Let $\Omega \subset M\setminus \partial M$ be a Lipschitz domain. Then for every $z\in \partial \Omega$ there exist an open neighborhood $U\subset M$ and a biLip\-schitz map $\psi\colon (0,1)\times [0,1)\to M$ such that $\psi((0,1)\times (0,1)) = U\cap \Omega$ and $\psi((0,1)\times\{0\}) = U\cap \partial\Omega$. Let $u\in W^{1,2}(\Omega, X)$. Then for almost every $s\in (0,1)$ the map $t\mapsto u\circ\psi(s,t)$ has an absolutely continuous representative which we denote by the same expression. The trace of $u$ is defined by $$\trace(u)(\psi(s,0)):= \lim_{t\searrow 0} (u\circ\psi)(s,t)$$ for almost every $s\in(0,1)$. It is shown in \cite{KS93} that the trace is independent of the choice of the map $\psi$ and defines an element of $L^2(\partial \Omega, X)$, that is, for some and hence every $x\in X$ the function $z\mapsto d(x, \trace(u)(z))$ is in $L^2(\partial \Omega)$.

\section{Morrey's lemma on $\varepsilon$--conformal mappings}

The aim of this section is to prove Theorem~\ref{thm:Morrey-epsilon-conformality}. We first establish Theorem~\ref{thm:Morrey-epsilon-Reshetnyak-disc-intro} which we restate for the convenience of the reader.

\bt\label{thm:Morrey-epsilon-conformality-disc}
 Let $X$ be a complete metric space and let $u\in W^{1,2}(D,X)$. Then for every $\varepsilon>0$ there exists a diffeomorphism $\varphi\colon D\to D$ such that $$E_+^2(u\circ\varphi) \leq \Area_{\mu_i}(u) + \varepsilon.$$
Moreover, the map $\varphi$ extends to a diffeomorphism of $\overline{D}$ and is conformal in a neighbourhood of the boundary.
\et

\begin{proof}
Let $\varepsilon>0$. We first claim that there exists $\delta>0$ such that the norms $s_z$, defined for almost every $z\in D$ by $$s_z(h) = \sqrt{\apmd u_z(h)^2 + \delta^2\cdot |h|^2},$$ satisfy 
\begin{equation}\label{eq:int-sz-almost-area}
\int_{D}\jac_{\mu_i}(s_z)\,dz \leq \Area_{\mu_i}(u) + \varepsilon.
\end{equation}
Indeed, by \eqref{eq:jac-energy} we have  $$\jac_{\mu_i}(s_z) \leq \mathbf{I}_+^2(s_z) = \mathbf{I}_+^2(\apmd u_z) + \delta^2$$ for almost every $z\in D$ and the function on the right-hand side is in $L^1(D)$. Since $\jac_{\mu_i}(s_z)$ converges to $\jac_{\mu_i}(\apmd u_z)$ as $\delta \searrow 0$ for almost every $z\in D$ the claim follows from dominated convergence.

Next, by the absolute continuity of the integral, there exists $L\geq 1$ such that the measurable set $A$ defined by $$A:= \left\{z\in \C: \text{$|z|\leq 1 - L^{-1}$ and $\mathbf{I}_+^1(\apmd u_z)\leq L$}\right\}$$ satisfies 
\begin{equation}\label{eq:small-int-energy-A}
\int_{D\setminus A}\mathbf{I}_+^2(\apmd u_z)\,dz \leq \varepsilon.
\end{equation}
Let $\mu\colon\C\to\C$ be the function defined as follows. For $z\in\C\setminus A$ set $\mu(z):= 0$. For $z\in A$ denote by $E_z$ the ellipse of maximal area contained in $$S_z:= \{h\in \C: s_z(h)\leq 1\}$$ and let $\mu(z)$ be the Beltrami coefficient of an orientation preserving linear isomorphism which takes $E_z$ to a round Euclidean ball. (The Beltrami coefficient is independent of the choice of such a linear map.) This defines a measurable function $\mu$ with compact support in $D$. It is furthermore not difficult to see that $\|\mu\|_\infty<1$. Indeed, if $z\in A$ then by definition $\mathbf{I}_+^1(\apmd u_z)\leq L$ and the ellipse $E_z$ satisfies $E_z \subset S_z\subset \sqrt{2}\cdot E_z$ by John's theorem. From this it follows that the eccentricity of $E_z$ is bounded by $K:= (2L^2\delta^{-2}+2)^{\frac{1}{2}}$ and we deduce from \eqref{eq:Beltrami-to-distortion} that $$\|\mu\|_\infty \leq \frac{K-1}{K+1}=:k<1.$$

Now, approximating $\mu$ by smooth functions via convolution with standard mollifiers and applying Egoroff's theorem and absolute continuity of the integral we obtain the following. There exist a smooth function $\tilde{\mu}\colon \C\to\C$ with compact support in $D$ and a measurable set $B\subset D$ satisfying  
\begin{equation}\label{eq:small-int-energy-B}
\int_{D\setminus B} \mathbf{I}_+^2(\apmd u_z)\,dz \leq \frac{\varepsilon}{K}
\end{equation}
 and such that $|\tilde{\mu}(z)|\leq k$ for every $z\in\C$ as well as 
\begin{equation}\label{eq:approx-unif-Beltrami-coeff}
|\mu(z) - \tilde{\mu}(z)|\leq (1-k^2)\cdot \frac{\varepsilon}{2+\varepsilon}
\end{equation} 
for every $z\in B$.

By Theorem~\ref{thm:Sol-Beltrami-qc-diffeo} there exists a diffeomorphism $\rho\colon\C\to\C$ which is $K$--quasicon\-formal and solves the Beltrami
equation $\rho_{\overline{z}} = \tilde{\mu} \rho_z$. Then the inverse $\rho^{-1}$ of $\rho$ is also $K$--quasiconformal, see Section~\ref{sec:Beltrami}.
Denote by $\varphi\colon \Omega\to D$ the restriction of $\rho^{-1}$ to $\Omega:= \rho(D)$. 
We will now estimate from above the energy of the composition $u\circ\varphi$. For this we first establish pointwise inequalities. Let $z\in D$ be such that the approximate metric derivative of $u$ at $z$ exists and set $w:= \rho(z)$. We distinguish cases as follows. Suppose first that $z\in B\cap A$ and let $T$ be an orientation preserving linear isomorphism which maps the ellipse $E_z$ of maximal area in $S_z$ to a round Euclidean ball. Then \eqref{eq:approx-unif-Beltrami-coeff} together with the composition formula \eqref{eq:composition-Beltrami-coeff} yields $$|\mu_{T\circ \varphi}(w)| = \left|\frac{\mu(z) - \tilde{\mu}(z)}{1- \mu(z)\overline{\tilde{\mu}(z)}}\right| \leq \frac{\varepsilon}{2+\varepsilon}$$ and hence $$\|T\circ D\varphi(w)\|^2 \leq (1+\varepsilon)\cdot \det(T\circ D\varphi(w))$$ by \eqref{eq:Beltrami-to-distortion}. Since the norm $s_z\circ T^{-1}$ is isotropic in the sense of Definition~\ref{def:isotropic-seminorm} we have $$\mathbf{I}_+^2(s_z\circ T^{-1}) = \jac_{\mu_i}(s_z\circ T^{-1})$$ and thus, together with the chain rule for the approximate metric derivative, we obtain that 
\begin{equation*}
 \begin{split}
  \mathbf{I}_+^2(\apmd (u\circ\varphi)_w) &\leq \mathbf{I}_+^2(\apmd u_z \circ T^{-1})\cdot \|T\circ D\varphi(w)\|^2\\
  &\leq (1+\varepsilon) \cdot \mathbf{I}_+^2(s_z\circ T^{-1}) \cdot \det(T\circ D\varphi(w))\\
  &= (1+\varepsilon)\cdot \jac_{\mu_i}(s_{\varphi(w)}) \cdot \det(D\varphi(w)).
 \end{split}
\end{equation*}
Next, we consider the case that $z\in B\setminus A$. Then $\mu(z)=0$ and thus \eqref{eq:approx-unif-Beltrami-coeff} and \eqref{eq:composition-Beltrami-coeff} imply $$\|D\varphi(w)\|^2\leq  (1+\varepsilon)\cdot \det(D\varphi(w))$$ and thus $$\mathbf{I}_+^2(\apmd (u\circ\varphi)_w) \leq (1+\varepsilon)\cdot \mathbf{I}_+^2(\apmd u_{\varphi(w)})\cdot \det(D\varphi(w)).$$ Finally, if $z\not\in B$ then the $K$--quasiconformality of $\varphi$ yields $$\mathbf{I}_+^2(\apmd (u\circ\varphi)_w) \leq K \cdot \mathbf{I}_+^2(\apmd u_{\varphi(w)})\cdot \det(D\varphi(w)).$$ From these pointwise inequalities we obtain by integrating and using the change of variables formula that
\begin{equation*}
 \begin{split}
 E_+^2(u\circ\varphi) &=\int_\Omega \mathbf{I}_+^2(\apmd (u\circ\varphi)_w)\,dw\\
 &\leq (1+\varepsilon) \int_{B\cap A} \jac_{\mu_i}(s_z)\,dz + (1+\varepsilon) \int_{B\setminus A} \mathbf{I}_+^2(\apmd u_z)\,dz\\
 &\quad + K\int_{D\setminus B} \mathbf{I}_+^2(\apmd u_z)\,dz\\
 &\leq (1+\varepsilon)\cdot [\Area_{\mu_i}(u) + \varepsilon] + (1+\varepsilon)\cdot \varepsilon + \varepsilon,
 \end{split}
\end{equation*}
where we have used \eqref{eq:int-sz-almost-area}, \eqref{eq:small-int-energy-A} and \eqref{eq:small-int-energy-B} in the last inequality. 

In order to complete the proof of the theorem we first observe that, after possibly precomposing with a conformal diffeomorphism from $D$ to $\Omega$, we may assume that $\Omega$ is the unit disc $D$. Notice that the Reshetnyak energy $E_+^2$ is conformally invariant. Since $\partial\Omega = \rho(S^1)$ is a smooth Jordan curve, a conformal diffeomorphism from $D$ to $\Omega$ extends to a diffeomorphism from $\overline{D}$ to $\overline{\Omega}$. Consequently, $\varphi$ can be extended to a diffeomorphism on $\overline{D}$. Finally, by construction, the function $\tilde{\mu}$ has compact support in $D$ and hence $\rho$ is conformal in a neighbourhood of the boundary $S^1$. In particular, $\varphi$ is conformal in a neighbourhood of the boundary as well.
\end{proof}

We use the result above to establish our main theorem.

\begin{proof}[Proof of Theorem~\ref{thm:Morrey-epsilon-conformality}]
Fix a Riemannian metric $g_0$ on $M$ and let $\varepsilon>0$. There exist pairwise disjoint smooth charts $(U_k, \psi_k)$ in $M$, $k=1,\dots, n$, which are all conformal with respect to $g_0$ and satisfy the following properties. Each $\psi_k$ is biLip\-schitz with $\psi_k(U_k) = D$, the sets $U_k$ do not intersect the boundary $\partial M$, and $$E_+^2(u|_K, g_0) \leq \frac{\varepsilon}{2},$$ where we have set $K:= M\setminus \bigcup_{k=1}^n U_k$.

Let $k\in\{1,\dots, n\}$ and notice that $u\circ\psi_k^{-1}\in W^{1,2}(D, X)$. By Theorem~\ref{thm:Morrey-epsilon-conformality-disc} there exists a diffeomorphism $\varphi_k\colon D\to D$ with
$$E_+^2(u\circ\psi_k^{-1}\circ \varphi_k)\leq \Area_{\mu_i}(u\vert_{U_k})+\frac{\varepsilon}{2n}$$
and such that $\varphi_k$ extends to a diffeomorphism of the closed unit disc and $\varphi_k$ is conformal outside a compact subset $A_k$ of $D$. We define $\varrho_k:= \varphi_k^{-1}\circ\psi_k$.

Let $V:= M\setminus \bigcup_{k=1}^n \varrho_k^{-1}(A_k)$ and note that $K\subset V$. Since $K$ is compact and $V$ is open there exist finitely many additional smooth charts $(U_k, \varrho_k)$, $k=n+1,\dots, m$, which are conformal with respect to $g_0$, contained in $V$ and cover $K$. It follows that $$\{(U_k,\varrho_k): k=1,\dots, m\}$$ is a smooth atlas on $M$ and all transition maps are conformal. Hence the uniformization theorem for compact Riemann surfaces implies the existence of a smooth metric $g$ on $M$ such that every chart $(U_k, \varrho_k)$ is conformal with respect to $g$. Moreover, $g$ can be chosen such that $(M,g)$ has constant sectional curvature $-1$, $0$, or $1$ and the boundary of $M$ is geodesic.

Finally, we have for each $k=1,\dots, n$ that $$E_+^2(u|_{U_k}, g) = E_+^2(u\circ\varrho_k^{-1}) \leq \Area_{\mu^i}(u|_{U_k}) + \frac{\varepsilon}{2n}.$$ Moreover, the identity map from $(M,g)$ to $(M,g_0)$ is conformal on $V$, thus $$E_+^2(u|_K, g) = E_+^2(u|_K, g_0) \leq \frac{\varepsilon}{2}.$$ Putting this together we see that $$E_+^2(u, g) = E_+^2(u|_K, g) + \sum_{k=1}^n E_+^2(u|_{U_k}, g) \leq \Area_{\mu^i}(u) + \varepsilon,$$
which finishes the proof.
\end{proof}

\section{Area minimizing surfaces and Courant's condition of cohesion}\label{sec:area-min-cond-cohesion}

In this section, we recall Courant's condition of cohesion and prove Theorem~\ref{thm:plateau-douglas-cond-cohesion}. We furthermore establish a weak analog of the theorem for the parametrized Hausdorff area.

Let $X$ be a complete metric space and let $\Gamma\subset X$ be the disjoint union of $k\geq 1$ rectifiable Jordan curves in $X$. Let $M$ be a smooth, compact, connected, orientable $2$--dimensional manifold with $k$ boundary components. We set $$a_{\mu^i}(M, \Gamma, X):= \inf\{\Area_{\mu^i}(v): v\in\Lambda(M, \Gamma, X)\}$$ and $$e(M,\Gamma, X):= \inf\{E_+^2(v, h): \text{ $v\in\Lambda(M, \Gamma, X)$ and $h$ Riemannian metric on $M$}\}.$$
The following result generalizes \cite[Theorem 1.1]{LW17-en-area}.

\bp\label{prop:area-min-energy}
We have 
\begin{equation}\label{eq:energy=area-inf}
 e(M,\Gamma, X) = a_{\mu^i}(M,\Gamma, X).
\end{equation}
 In particular, if $u\in \Lambda(M, \Gamma, X)$ and a Riemannian metric $g$ on $M$ are such that $$E_+^2(u,g) = e(M,\Gamma, X)$$ then $u$ satisfies $$\Area_{\mu^i}(u) =  a_{\mu^i}(M,\Gamma, X)$$  and $u$ is infinitesimally isotropic with respect to $g$.
\ep

\begin{proof}
 The equality \eqref{eq:energy=area-inf} is a consequence of Theorem~\ref{thm:Morrey-epsilon-conformality} and \eqref{eq:area-energy-ineq}. Suppose now that $u\in \Lambda(M, \Gamma, X)$ and a Riemannian metric $g$ on $M$ are such that $$E_+^2(u,g) = e(M,\Gamma, X).$$
 Corollary~\ref{cor:energy-min-inf-isotropic} implies that $u$ is infinitesimally isotropic with respect to $g$ and hence $$\Area_{\mu^i}(u) = E_+^2(u,g) = e(M,\Gamma, X) = a_{\mu^i}(M,\Gamma, X)$$ by Remark~\ref{rem:area-engery-equality-inf-isotrop} and \eqref{eq:energy=area-inf}. 
\end{proof}

\bd
Given $\eta>0$, a map $u\colon M\to X$ is called $\eta$--cohesive if $u$ is continuous and for every non-contractible closed curve $c$ in $M$ the length of $u\circ c$ is no smaller than $\eta$.
A family of maps from $M$ to $X$ is said to satisfy the condition of cohesion if there exists $\eta>0$ such that each element of the family is $\eta$--cohesive.
\ed

The condition of cohesion was introduced by Courant \cite{Cou40} and used in \cite{Shi39} and \cite{TT88}. It is for example satisfied when the maps are incompressible in the sense of of Schoen-Yau \cite{SY79}.

\begin{proof}[Proof of Theorem~\ref{thm:plateau-douglas-cond-cohesion}]
Let $(u_n)$ be a $\mu^i$--area minimizing sequence in $\Lambda(M,\Gamma, X)$ which satisfies the condition of cohesion. By Theorem~\ref{thm:Morrey-epsilon-conformality} there exists, for each $n\in\N$, a Riemannian metric $g_n$ on $M$ such that $$E_+^2(u_n,g_n)\leq \Area_{\mu^i}(u_n) + \frac{1}{n},$$ so the sequence $(u_n)$ is also an $E_+^2$--energy minimizing sequence in $\Lambda(M,\Gamma, X)$ by \eqref{eq:energy=area-inf}. Since $(u_n)$ satisfies the condition of cohesion, \cite[Theorem 8.2]{FW-Plateau-Douglas} implies the existence of some $u\in\Lambda(M,\Gamma, X)$ and of a Riemannian metric $g$ such that $$E_+^2(u,g) = e(M,\Gamma, X)$$ and $u$ is infinitesimally isotropic with respect to $g$. Moreover, the Riemannian metric $g$ can be chosen in such a way that $(M,g)$ has constant curvature $-1$, $0$, $1$ and that $\partial M$ is geodesic. Proposition~\ref{prop:area-min-energy} shows that $u$ is $\mu^i$--area minimizing in $\Lambda(M,\Gamma, X)$. The proof is complete.
\end{proof}

We can furthermore provide a weak analog of Theorem~\ref{thm:plateau-douglas-cond-cohesion} when the intrinsic Riemannian area is replaced by the Hausdorff area. Recall that the {\it Hausdorff area} of a map $u\in W^{1,2}(M,X)$ is given by $$\Area(u)= \int_M \jac(\apmd u_z)\,d\hm^2_g(z),$$ where the (Hausdorff) jacobian of a semi-norm $s$ on Euclidean $\R^2$ is defined by $$\jac(s):=\left\{\begin{array}{cl} 
 \hm^2_{(\R^2,s)}(e_1\wedge e_2) & \text{if $s$ is norm}\\
 0&\text{otherwise.}
\end{array}\right.$$
As before we identify the tangent space $V=(T_zM, g(z))$ with Euclidean $\R^2$ via a linear isometry to define the jacobian of a semi-norm on $V$. If $u$ satisfies Lusin's property (N) then the area formula \cite{Kir94}, \cite{Kar07} for metric space valued maps yields $$\Area(u) = \int_X\#\{z\in \Omega: u(z) = x\}\,d\hm^2_X(x).$$ We have the inequality 
\begin{equation}\label{eq:Hausdorff-intrinsic-Riem-comparison}
\frac{\pi}{4}\cdot \Area_{\mu^i}(u) \leq \Area(u) \leq \Area_{\mu^i}(u)
\end{equation}
 for every $u\in W^{1,2}(M, X)$, see \cite[Section 2.4]{LW17-en-area}.
While the Hausdorff area is different from the intrinsic Riemannian area in the generality of normed or metric spaces the two notions agree for example when $X$ is a Riemannian manifold, or more generally, a metric space with property (ET) of `Euclidean tangents' in the sense of \cite{LW15-Plateau}.

\bp\label{prop:Plateau-Douglas-cohesion-Hausdorff-weak}
 Suppose the metric space $X$ is proper and $\Lambda(M, \Gamma, X)$ is not empty. Suppose furthermore that there is an area minimizing sequence (with respect to the Hausdorff area) in $\Lambda(M, \Gamma, X)$ satisfying the condition of cohesion. Then there exist $u\in\Lambda(M, \Gamma, X)$ such that $$\Area(u) =  \inf\{\Area(v): v\in\Lambda(M, \Gamma, X)\}.$$
\ep

Unlike in Theorem~\ref{thm:plateau-douglas-cond-cohesion} we do not obtain a good parametrization. This is in contrast to the situation in our paper \cite{FW-Plateau-Douglas}, where we solved the Plateau-Douglas problem under the Douglas condition. Notice that minimizers with respect to the Reshetnyak energy need not minimize the Hausdorff area in general, see \cite[Proposition 11.6]{LW15-Plateau}. At present, we do not know whether there exists a definition of energy in the sense of \cite{LW17-en-area} for which energy minimizers are Hausdorff area minimizers.

\begin{proof}
 Let $(u_n)$ be a Hausdorff area minimizing sequence in $\Lambda(M,\Gamma, X)$ which satisfies the condition of cohesion. By Theorem~\ref{thm:Morrey-epsilon-conformality} there exists, for each $n\in\N$, a Riemannian metric $g_n$ on $M$ such that $$E_+^2(u_n,g_n)\leq \Area_{\mu^i}(u_n) + 1.$$ It follows from \eqref{eq:Hausdorff-intrinsic-Riem-comparison} that the energies $E_+^2(u_n, g_n)$ are uniformly bounded. Arguing exactly as in the proof of \cite[Theorem 8.2]{FW-Plateau-Douglas} one shows the existence of $u\in\Lambda(M,\Gamma, X)$ such that a subsequence of the sequence $(u_n)$ converges to $u$ in $L^2(M, X)$ with respect to some fixed Riemannian metric on $M$. By the lower semi-continuity of the Hausdorff area \cite[Corollary 5.8]{LW15-Plateau} it follows that $u$ is an area minimizer in $\Lambda(M,\Gamma, X)$.
\end{proof}

\def\cprime{$'$} \def\cprime{$'$} \def\cprime{$'$}

\end{document}